\documentclass[11pt]{amsart}

\usepackage{amsfonts,amsmath}
\usepackage{amssymb, latexsym}
\usepackage{amscd,amsthm}
\input xy
\xyoption{all}

\newcommand{\marginlabel}[1]%
  {\mbox{}\marginpar{\raggedleft\hspace{0pt}\bfseries\sf#1}}

\def\CC{{\mathbb C}}

\def\OO{{\mathcal O}}
\def\R{\mathbf{R}}

\def\D{\mathbf{D}}

\def\F{\mathcal{F}}

\def\G{\mathcal{G}}

\def\I{\mathcal{I}}

\def\Pic0{{\rm Pic}^0(X)}

\theoremstyle{plain}

\newtheorem{theorem}{Theorem}[section]
\newtheorem{theoremalpha}{Theorem}

\newtheorem{conjecturealpha}[theoremalpha]{Conjecture}

\newtheorem{proposition/example}[theorem]{Proposition/Example}
\newtheorem{proposition}[theorem]{Proposition}
\newtheorem{corollary}[theorem]{Corollary}
\newtheorem{lemma}[theorem]{Lemma}

\theoremstyle{definition}
\newtheorem{definition}[theorem]{Definition}
\newtheorem{remark}[theorem]{Remark}
\newtheorem{example}[theorem]{Example}

\newtheorem{conjecture/question}[theorem]{Conjecture/Question}

\newtheorem{remark/definition}[theorem]{Remark/Definition}
\newtheorem{definition/notation}[theorem]{Definition/Notation}

\theoremstyle{remark}

\begin{document}

\title{Generic vanishing and minimal cohomology classes on abelian varieties}

\author[G. Pareschi]{Giuseppe Pareschi}
\address{Dipartamento di Matematica, Universit\`a di Roma, Tor Vergata, V.le della
Ricerca Scientifica, I-00133 Roma, Italy}
\email{{\tt pareschi@mat.uniroma2.it}}

\author[M. Popa]{Mihnea Popa}
\address{Department of Mathematics, University of Chicago,
5734 S. University Ave., Chicago, IL 60637, USA }
\email{{\tt mpopa@math.uchicago.edu}}

\thanks{2000\,\emph{Mathematics Subject Classification:} 14K12, 14F17}
\keywords{Principally polarized abelian
varieties, minimal classes, generic vanishing}
\thanks{MP was partially supported by the NSF grant DMS 0500985
and by an AMS Centennial Fellowship.}


\begin{abstract}
We establish a -- and conjecture further -- relationship between the
existence of subvarieties representing minimal cohomology classes on
principally polarized abelian varieties, and the generic vanishing
of the cohomology of twisted ideal sheaves. The main ingredient is the
Generic Vanishing criterion established in \cite{pp3}, based on the
Fourier-Mukai transform.
\end{abstract}
\maketitle


\setlength{\parskip}{.1 in}

\markboth{G. PARESCHI and M. POPA}
{\bf Generic vanishing and minimal cohomology classes}

\section{Introduction}

This paper is concerned with a relationship between the existence of subvarieties of principally polarized abelian varieties (ppav's)  having minimal cohomology class and the (generic) vanishing of certain sheaf cohomology, based on the Generic Vanishing criterion studied in \cite{pp3}. This is in analogy with the well-known equivalence between a subvariety in projective space being of minimal degree and its ideal sheaf being Castelnuovo-Mumford $2$-regular (cf. e.g. \cite{eg}). All the statements below are over $\CC$. We first state the following conjecture.

\begin{conjecturealpha}\label{a}
Let $(A, \Theta)$ be a indecomposable ppav of dimension $g$, and let $X$ be a geometrically nondegenerate closed reduced subscheme of $A$ of pure dimension $d\le g-2$.
The following are equivalent:
\begin{enumerate}
\item $X$ has minimal cohomology class, i.e. $[X] = \frac{\theta^{g-d}}{(g-d)!}$.
\item $\I_X (\Theta)$ is a $GV$-sheaf.
\item $\I_X(2\Theta)$ satisfies $IT_0$.
\item $\OO_X(\Theta)$ is $M$-regular, and $\chi(\OO_X(\Theta))=1$.
\item Either $(A, \Theta)$ is the polarized Jacobian of a smooth projective curve of genus $g$ and
$X$ is $+$ or $-$ an Abel-Jacobi embedded copy of $W_d( C)$, or $g= 5$, $d=2$, $(A, \Theta)$ is the
intermediate Jacobian of a smooth cubic threefold and $X$ is $+$ or $-$ a translate of the Fano surface of lines.
\end{enumerate}
\end{conjecturealpha}

The condition that $\I_X(\Theta)$ be a $GV$-sheaf means that it satisfies a condition analogous to the
Green-Lazarsfeld result for canonical bundles \cite{gl}, namely
$${\rm codim}~\{ P_{\alpha} \in {\rm Pic}^0 (A) ~| ~ h^i (\I_X(\Theta) \otimes P_{\alpha}) \neq 0\} \ge i
{\rm ~for~all~} i.$$
(Cf. \S2 for the notions used in (3) and (4).) The condition of being geometrically nondegenerate is the weakest nondegeneracy condition
one usually considers. It is defined in \cite{ran} \S II, together with the stronger condition of being nondegenerate (cf. \S4) -- note however that subvarieties of minimal class are already known to satisfy the stronger version (cf. \cite{ran}, Corollaries II.2 and II.3), so anything giving (1) implies a posteriori that $X$ is in fact \emph{nondegenerate}.

Part (5), more precisely its equivalence with (1), is of course not
our conjecture, but rather was formulated in low dimensions by Beauville \cite{beauville1} and 
Ran \cite{ran}, and then in general by Debarre \cite{debarre1}. 
In \S2 we explain that most of the implications between the cohomological
statements -- namely (2) $\iff$ (4) $\Rightarrow$ (3) -- are already
known or follow quickly from the definitions. We don't know how to
prove that (3) implies the other two. We have already conjectured
that (3) is equivalent to (1) before (cf. \cite{pp2}, where the
analogy with the Castelnuovo-Mumford regularity picture is
suggested); we will see here however that condition (2) is more
natural in connection with minimal classes. Note also that in
\cite{pp1} Proposition 4.4 it is shown that (5) $\Rightarrow$ (4)
for $W_d$'s in Jacobians, while in \cite{hoering} the same thing is
proved for the Fano surface of a smooth cubic threefold. Thus (5) is
known to imply all (1) -- (4). Note also that in dimension four we
do know the equivalence of (1), (2), (4) and (5), and the fact that
they imply (3): indeed, a result of Ran \cite{ran} asserts the
equivalence of (1) and (5) in Conjecture \ref{a}, while Theorem
\ref{b} below gives that (2) implies (1).

Our main results here are concerned with what is implied by the Generic Vanishing condition. The first is a proof of a slightly stronger (2) $\Rightarrow$ (1) in the Conjecture.\footnote{Note that indecomposability plays no role here, and in fact in anything that is not related to (5).}

\begin{theoremalpha}\label{b}
Let $X$ be a geometrically nondegenerate closed reduced subscheme of pure dimension $d$ of a ppav $(A, \Theta)$ of dimension $g$. If $\I_X (\Theta)$ is a $GV$-sheaf, then $X$ is Cohen-Macaulay and $[X] = \frac{\theta^{g-d}}{(g-d)!}$.
\end{theoremalpha}

The second is a sheaf cohomology criterion for detecting Jacobians. It is the implication (2) $\Rightarrow$ (5) in the Conjecture in the cases $d = 1$ and $d = g-2$. The case $d=1$ is in fact an
immediate consequence of Theorem \ref{b} and the Matsusaka-Ran criterion, so the main content is that the same applies for codimension two subvarieties.

\begin{theoremalpha}\label{c}
Let $X$ be a geometrically nondegenerate equidimensional reduced subscheme of a $g$-dimensional indecomposable ppav $(A, \Theta)$, of dimension either $1$ or $g-2$ respectively. If $\I_X (\Theta)$ is a $GV$-sheaf, then $(A, \Theta)$ is the polarized Jacobian of a smooth projective curve $C$ of genus $g$, and $X$ is $+$ or $-$ an Abel-Jacobi embedded copy of either $C$ or $W_{g-2}(C)$ respectively.
\end{theoremalpha}

There are two new tools that are used for the proofs of Theorems
\ref{b} and \ref{c}. The main one is \cite{pp3} Theorem F, which
equates the $GV$-condition for an object $\F$ in the derived
category $\D(A)$ with the vanishing of suitable components of the
Fourier-Mukai transform of the Grothendieck dual of $\F$ (cf. \S2
below for a review). The other is to consider systematically the
locus of theta-translates containing $X$, which we call the
theta-dual $V(X)$ of $X$ in $\widehat A$.
 As we will show in \S5, such locus
supports sheaves which are components of the Fourier-Mukai transform
of naturally defined complexes. Roughly speaking, we relate the
$GV$-condition for $\I_X(\Theta)$ to the cohomology class of $V(X)$
via Grothendieck-Riemann-Roch.

In \S8 we formulate natural geometric conditions on $V(X)$ which should be equivalent to those in Conjecture \ref{a}. We also note that a proof of (1) $\Rightarrow$ (2) would have a consequence for
the implication (1) $\Rightarrow$ (5) as well: one would need to check this only for $X$ of dimension $d \le [\frac{g}{2}]$.

\noindent {\bf Acknowledgements.} We are grateful to O. Debarre for
pointing out how to arrange things in the proof so that we could
assume only geometric nondegeneracy. We also thank Ch. Hacon for
suggesting Example \ref{counterexample}, and to A. H\"oring for
showing us a preliminary version of \cite{hoering}. Finally, we thank the 
referee for numerous comments which improved the math and the exposition.

\section{Preliminaries}

We recall the main Fourier-Mukai terminology and results used in the sequel.
If $A$ is an abelian variety of dimension $g$ and $\widehat{A}$ is the dual abelian variety, let
$\mathcal{P}$ be a normalized Poincar\'e line bundle on $A\times \widehat{A}$. The Fourier-Mukai
functor \cite{mukai1} is
$$\mathbf{R}\widehat{\mathcal{S}}:{\bf D}(A )\rightarrow {\bf D}(\widehat A),$$
the derived functor induced by $\widehat{\mathcal{S}}: {\rm Coh(A)}\rightarrow {\rm Coh}(\widehat{A})$,
where $\widehat{\mathcal{S}} (\F) = {p_{\widehat A}}_* (p_A^* \F \otimes \mathcal{P})$.
We also consider $\mathbf{R}\mathcal{S}:{\bf D}(\widehat A )\rightarrow {\bf D}(A)$,
defined analogously. Mukai's main result (cf. \cite{mukai1}, Theorem 2.2) is that $\R \widehat{\mathcal{S}}$ is an equivalence of derived categories and
\begin{equation}\label{inversion}
\mathbf{R}\mathcal{S}\circ\mathbf{R}\widehat{\mathcal{S}}\cong (-1_A)^{*}
[-g] ~{\rm ~and~}~ \mathbf{R}\widehat{\mathcal{S}}\circ\mathbf{R}\mathcal{S}\cong (-1_{\widehat{A}})^{*} [-g].
\end{equation}

An object $\F$ in $\D (A)$ is said to satisfy $WIT_i$ (Weak Index Theorem with index $i$) if
$R^j\widehat{\mathcal{S}}(\F)=0, {\rm ~for~ all~} j\neq i $. In this case
$R^{i}\hat{\mathcal{S}}(\F)$ is denoted by $\widehat{\F}$ and called the $\it{Fourier~
transform}$ of $\F$. Note that $\mathbf{R}\widehat{\mathcal{S}}(\F)\cong \widehat{\F}[-i]$.
Moreover, $\F$ satisfies the stronger $IT_i$ (Index Theorem with index $i$) if
$$H^{j}(A, \F\otimes P_{\alpha})=0 {\rm~ for~ all~} \alpha\in \widehat A  {\rm~ and~ all~} j\neq i,$$
where (as in the rest of the paper) we denote by $P_{\alpha} \in
{\rm Pic}^0 (A)$ the line bundle corresponding to the point
$\alpha\in \widehat A$. Reversing the role of $A$ and $\widehat A$,
and using the functor $\R\mathcal{S}$ instead, all the previous
notions/notation can be defined for objects in ${\bf D}(\widehat
A)$.

Given an object $\F$ in $\D(A)$  we use the notation
$$\R \Delta \F : = \R \mathcal{H}om (\F , \OO_A),$$
and similarly on $\widehat A$. Note that the Grothendieck dualizing
functor applied to $\F$ is a shift of this, namely $\R \Delta \F
[g]$. Grothendieck duality, applied to the present context, states
that  (cf. \cite{mukai1} (3.8)):
\begin{equation}\label{gd}
\R\Delta \circ \R \widehat{\mathcal{S}} \cong ((-1_{\widehat A})^* \circ \R
\widehat{\mathcal{S}} \circ \R \Delta) [g], \qquad \R\Delta \circ \R  \mathcal{S}
\cong ((-1_{ A})^* \circ \R  \mathcal{S} \circ \R \Delta) [g].
\end{equation}

An object $\F$ in $\D(A)$ is called a \emph{$GV$-object} (cf. \cite{pp3}\footnote{In
\cite{pp3} we define $GV$-objects with respect to a given Fourier-Mukai functor. Since in the present paper we will use only the functor $\R \widehat{\mathcal{S}}$, we suppress this from the notation.}) if 
$${\rm codim}~{\rm Supp}(R^i \widehat{\mathcal{S}} (\F) )\geq i {\rm ~for~all~} i.$$
More generally, for any integer $k \ge 0$, an object $\F$ in $\D(A)$ is called a \emph{$GV_k$-object} if  
$${\rm codim}~{\rm Supp}(R^i \widehat{\mathcal{S}} ( \F))\geq i- k {\rm ~for~all~} i.$$
(So $GV = GV_0$.). Although not strictly necessary for the purpose of this paper, it is worth recalling that, by \cite{pp3} Lemma 3.8, a coherent sheaf $\F$ is $GV_k$ if and only if the more familiar
 Green-Lazarsfeld condition
$${\rm codim}_{\widehat A} V^i (\F)\geq i-k$$
is satisfied for all $i$, where $V^i (\F) : = \{\alpha \in \widehat
A ~| ~ H^i( A,\F \otimes P_{\alpha}) \neq 0\}$ is the $i$-th
cohomological support locus of $\F$. (Same for any object in $\D(A)$
if we use hypercohomology ${\bf H}^i$.) The main technical tool of the paper
is a particular case of the Generic Vanishing criterion \cite{pp3}
Theorem F. This holds in a much more general context --  for
sake of self-containedness, here we state the particular result
needed in the paper, for which we provide a short {\it ad hoc} proof.

\begin{theorem}\label{gv}
With the notation above, let $\F$ be an object in $\D(A)$, with cohomologies only in non-negative degrees. The following are equivalent:
\newline
\noindent
(a) $\F$ is a $GV_k$-object.
\newline
\noindent 
(b) ${\bf H}^i (A, \F \otimes \widehat{L^{-1}}) = 0$ for
$i \notin [0,k]$, for any sufficiently positive ample line bundle
$L$ on $\widehat A$.
\newline
\noindent 
(c)  $R^i \widehat{\mathcal{S}} ( \R \Delta \F) = 0$ for all $i\notin [g-k, g]$.
\end{theorem}

\begin{proof} 
(b) $\Leftrightarrow$ (c). To begin with, it is a standard fact that, for a fixed $i$, the vanishing
of $R^i \widehat{\mathcal{S}} ( \R \Delta \F)$ is equivalent to
\begin{equation}\label{intermediate}{\bf
H}^i (\widehat A, \R\widehat{\mathcal{S}}(\R \Delta \F )\otimes L) = 0
\quad\hbox{for any sufficiently positive ample line bundle $L$ on
$\widehat A$.}\footnote{The proof in brief: since $L$ is
sufficiently positive, by Serre vanishing the hypercohomology
spectral sequence degenerates, hence ${\bf H}^i (\widehat A, \R\widehat{\mathcal{S}}(\R
\Delta \F )\otimes L)\cong H^0(\widehat A, R^i\widehat{\mathcal{S}}(\R \Delta \F )\otimes
L)$. But, again by Serre's theorem, the right hand side vanishes if and
only if $R^i\widehat{\mathcal{S}}(\R \Delta \F )$ does.}
\end{equation} 
We rewrite
$${\bf H}^i (\widehat A, \R\widehat{\mathcal{S}}(\R \Delta \F )\otimes L)={\rm 
Ext}^i_{\D(\widehat A)} (L^{-1},\R\widehat{\mathcal{S}}(\R \Delta \F ))={\rm
Hom}_{\D(\widehat A)}(L^{-1},\R\widehat{\mathcal{S}}(\R \Delta \F) [i]).$$ 
Since $\R\mathcal{S}$ is an equivalence, the last ${\rm Hom}$ 
is isomorphic to ${\rm Hom}_{\D(A)}(\R\mathcal{S}
(L^{-1}), \R\mathcal{S} (\R\widehat{\mathcal{S}}(\R \Delta \F ))[i])$. Now 
$L^{-1}$ is $IT_g$ and hence $\R\widehat{\mathcal{S}}
(L^{-1}) = \widehat{L^{-1}}[-g]$. Moreover, by Mukai's theorem
(\ref{inversion}) we have that $\R\mathcal{S} (\R\widehat{\mathcal{S}}(\R \Delta \F ))\cong
(-1)^*_A\R\Delta\F[-g]$. Therefore
$${\rm Hom}_{\D(A)}(\R\mathcal{S}
(L^{-1}),\R\mathcal{S} (\R\widehat{\mathcal{S}}(\R \Delta \F))[i])\cong {\rm
Hom}_{\D(A)}(\widehat{L^{-1}},(-1)^*_A\R\Delta\F[i]).$$ 
The right hand side is isomorphic to ${\rm 
Ext}^i(\widehat{L^{-1}},(-1)^*_A\R\Delta\F)$, which by Grothendieck-Serre
duality is isomorphic to ${\bf H}^{g-i}(A,\widehat{L^{-1}}\otimes(-1)^*_A\F)$. 
In conclusion, we have proved that $R^i \widehat{\mathcal{S}} ( \R \Delta \F)$ vanishes if and
only if the hypercohomology group ${\bf
H}^{g-i}(A,\widehat{L^{-1}}\otimes(-1)^*_A\F)$ does, for any
sufficiently positive ample line bundle $L$ on $\widehat A$. This
last condition is clearly equivalent to the vanishing of ${\bf
H}^{g-i}(A,\widehat{L^{-1}}\otimes\F)$, which proves (b) $\Leftrightarrow$ (c).

\noindent (a) $\Rightarrow$ (b). Since $L^{-1}$ is $IT_g$, the
transform $\widehat{L^{-1}}$ is locally free. The required vanishing for $i <0$ is obvious since 
$\F$ has cohomologies only in non-negative degrees.
Now by Grothendieck duality (\ref{gd}) we have 
$(\widehat{L^{-1}})^\vee\cong  (-1)_{A}^* \widehat L$. Therefore
$${\bf H}^{i}(A,\widehat{L^{-1}}\otimes\F)\cong{\rm Ext}_{\D(A)}^i ((-1)^*_A\widehat
L, \F)\cong{\rm Hom}_{\D(A)}((-1)^*_A\widehat L, \F[i]).$$ 
Since $\R\widehat{\mathcal{S}}$ is an equivalence, the last ${\rm
Hom}$ is isomorphic to ${\rm Hom}_{\D(\widehat
A)}(\R\widehat{\mathcal{S}} ((-1)^*_A\widehat L),\R\widehat{\mathcal{S}} (\F)[i]),$
which is in turn isomorphic to
 $${\rm Hom}_{\D(\widehat
A)}( L[-g],\R\widehat{\mathcal{S}} (\F)[i])\cong {\rm Ext}^{g+i}_{D(\widehat
A)}(L,\R\widehat{\mathcal{S}} (\F))\cong {\bf H}^{g+i}(\widehat A, L^{-1}\otimes
\R\widehat{\mathcal{S}} (\F)).$$ 
Here we used Mukai inversion
(\ref{inversion}) to deduce that $\R\widehat{\mathcal{S}} ((-1)^*_A\widehat
L) = (\R\widehat{\mathcal{S}}\circ (-1)^*_A\circ \R\mathcal{S})( L) \cong L[-g]$.
Therefore we are reduced to proving that ${\bf H}^{g+i}(\widehat A, L^{-1}\otimes
\R\widehat{\mathcal{S}} (\F))=0$ for $i > k$ as soon as $\F$ is a
$GV_k$-object. This follows easily from the hypercohomology spectral
sequence $$E^{ij}_2 := H^j(\widehat A, L^{-1}\otimes
R^i\widehat{\mathcal{S}}(\F))\Rightarrow {\bf H}^{i+j}(\widehat A, L^{-1}\otimes
\R \widehat{\mathcal{S}}(\F)),$$ since $GV_k$ means that the
$R^i\widehat{\mathcal{S}}(\F)$'s are supported in dimension $\le g-i + k$, for any $i$.

\noindent (c) $\Rightarrow$ (a). Since $\R\Delta$ is an involution
on $\D(A)$, Grothendieck duality (\ref{gd}), applied to
$\R\Delta\F$, gives 
$$\R\Delta (\R \widehat{\mathcal{S}}(\R\Delta\F))
\cong (-1_{\widehat A})^*  \R \widehat{\mathcal{S}}(\F) [g].$$ 
This gives rise to a spectral sequence
$$E^{ij}_2: = {\mathcal E}xt^{i+j}(R^j\widehat{\mathcal{S}} (\R\Delta
\F [g]), \OO_{\widehat A})\Rightarrow (-1)_{\widehat A}^*R^i\widehat{\mathcal{S}}(\F).$$ 
The hypothesis is that $R^j\widehat{\mathcal{S}}(\R\Delta \F)$ vanish
for $j\not\in [g-k,g]$, i.e. that the $R^j:=R^j\widehat{\mathcal{S}} (\R\Delta
\F[g])$ vanish for $j\not\in [-k,0]$. One knows that
$${\rm codim} ~{\rm Supp}(\mathcal{E}xt^{i + j}(R^j, \OO_Y)) \ge i + j$$
for all $i$ and $j$. This is a particular case of the following general fact: 
for any $h$, the ${\mathcal E}xt^h(\G,\OO_X)$-sheaf associated to a sheaf $\G$ 
on a smooth variety $X$ is always supported in codimension $\ge h$ (see e.g.
\cite{oss}, Ch.II, Lemma 1.1.2).  Since the only non-zero
$R^j$-sheaves are for $j \ge -k$, we have that the codimension of
the support of  every $E_{\infty}$ term of the above spectral
sequence is at least $i -k$. This implies immediately that
$R^i\widehat{\mathcal{S}} (\F) =0$ for $i<0$ and ${\rm codim}~{\rm Supp}(R^i
\widehat{\mathcal{S}}(\F)) \geq i-k$, for all $i\ge 0$.
\end{proof}

An important particular case is that of $k = 0$, when the result says that $\F$ is $GV$ if and only if $\R \Delta \F$ satisfies $WIT_g$. In this case, the fact that (c) is equivalent to (b) and implies (a) appears already in \cite{hacon}.

Finally, we recall that a sheaf on $A$ is called \emph{M-regular}
(cf. \cite{pp1}) if $${\rm codim}~{\rm Supp}(R^i \widehat{\mathcal{S}} (\F) )> i {\rm ~for~all~} i> 0.$$

\section{Implications between the cohomological statements}

The implication (2) $\Rightarrow$ (3) in Conjecture \ref{a} follows
from the general Lemma below, which is proved exactly as the
implication (1) $\Rightarrow$ (2) of Theorem \ref{gv} (note that 
$\widehat{\OO_{\widehat A}(-\Theta)} = (-1_A)^* \OO_A(\Theta)$ -- cf. \cite{mukai1}, Proposition 3.11). In \S7 below we will see some extra implications of the $GV$-hypothesis when $\F$
is (the twist of) an ideal sheaf.

\begin{lemma}\label{arbitrary}
Let $\F$ be a $GV$-sheaf on a ppav $(A, \Theta)$. Then $\F (\Theta)$ satisfies $IT_0$.
\end{lemma}

\begin{example}\label{counterexample}
Given Lemma \ref{arbitrary}, it is natural to ask whether the implication (3) $\Rightarrow$ (2) might also hold for an arbitrary coherent sheaf $\F$ instead of $\I_X(\Theta)$. We give an example showing that this is not the case, so this implication (if true) should be more geometric.
We thank Christopher Hacon for suggesting this example. Consider $(A, \Theta)$ any ppav, say with $\Theta$ symmetric for simplicity, $k \ge 2$ an integer, and $\phi_k : A \rightarrow A$ the
map given by multiplication by $k$. Then $\F : = {\phi_k}_* \OO_A (-\Theta)$ is a sheaf which is not
$GV$, but $\F (\Theta)$ satisfies $IT_0$. Indeed $\F(\Theta)$ satisfies $IT_0$ if and only if
$\OO_A(-\Theta) \otimes \phi_k^* \OO_A(\Theta)$ does so, which is obviously true since $\phi_k^*
\OO_A(\Theta) \cong \bigoplus \OO_A (k^2 \Theta)$.  On the other hand, $\F$ is not $GV$ since
$$H^g (A, \F \otimes P_{\alpha}) \cong H^g (A, \OO_A(-\Theta) \otimes \phi_k^* P_{\alpha}),$$
which is non-zero for all $\alpha$.
\end{example}

The equivalence (2) $\iff$ (4) follows directly from the definitions, without any assumptions on $X$.

\begin{lemma}\label{regularity_equivalence}
Let  $(A, \Theta)$ be a ppav and $X$ a closed subscheme of $A$. The following are equivalent:

\noindent
(a) $\I_X(\Theta)$ is a $GV$-sheaf.

\noindent
(b) $\OO_X(\Theta)$ is $M$-regular and $\chi(\OO_X(\Theta))=1$.
\end{lemma}
\begin{proof}
Applying the Fourier-Mukai functor to the exact sequence
$$0 \longrightarrow \I_X(\Theta) \longrightarrow \OO_A(\Theta) \longrightarrow \OO_X(\Theta)
\longrightarrow 0,$$
we obtain that $R^i\widehat{\mathcal S}(\OO_X(\Theta)) \cong  R^{i +1}\widehat{\mathcal S}
(\I_X(\Theta))$ for all $i \ge 1$, and
$$0 \longrightarrow \OO_A(-\Theta) \longrightarrow R^0\widehat{\mathcal S}(\OO_X(\Theta)) \longrightarrow R^1\widehat{\mathcal S}(\I_X(\Theta)) \longrightarrow 0,$$
since $R^0\widehat{\mathcal S}(\I_X(\Theta))$ being torsion must therefore be zero.
Thus if $\I_X(\Theta)$ is a $GV$-sheaf, then $\OO_X(\Theta)$ is $M$-regular and $R^0\widehat{\mathcal S}(\OO_X(\Theta))$ has generic rank $1$, which by base change is equivalent to
 $h^0(A, \OO_X(\Theta)\otimes P_{\alpha})=1$ for $P_{\alpha} \in {\rm Pic}^0 (A)$ general, i.e. $\chi (\OO_X(\Theta)) = 1$. Conversely, if (b) holds than the same reasoning shows that we need to worry
 only about $R^1\widehat{\mathcal S}(\I_X(\Theta))$. But this must be supported on  a proper subset,
 since the argument above can be completely reversed.
\end{proof}

\section{The theta-dual of a subvariety of a ppav}

Let $(A, \Theta)$ be a ppav of dimension $g$. Consider a closed reduced subscheme $X$ of $A$, of dimension $d \le g-2$. We put a natural scheme structure on the locus of theta-translates containing
$X$.

\begin{lemma}\label{identification}
For any $\alpha \in\widehat A$ there is a canonical identification
$$(-1_{\widehat A} )^* R^g\widehat{\mathcal S}({\bf R}\Delta(\I_X(\Theta))) \otimes_{\OO_{\widehat{A} , \alpha}} \CC(\alpha)\cong
H^0 (A, \I_X(\Theta) \otimes P_{\alpha})^\vee.$$
\end{lemma}
\begin{proof}
We apply base change for complexes, as in \cite{ega3} 7.7:  since ${\bf H}^{g+1} (A, {\bf R}\Delta(\I_X(\Theta)) \otimes P_{\alpha}) = 0$ for all $P_{\alpha} \in {\rm Pic}^0 (A)$, we have for any $\alpha$ a natural isomorphism
$$(-1_{\widehat A} )^* R^g\widehat{\mathcal S}({\bf R}\Delta(\I_X(\Theta)))\otimes_{\OO_{\widehat{A} , \alpha}} \CC(\alpha)\cong
{\bf H}^g (A, {\bf R}\Delta(\I_X(\Theta)) \otimes P_{\alpha}^{-1}).$$
By Grothendieck-Serre duality, this hypercohomology group is isomorphic
to $H^0 (A, \I_X(\Theta) \otimes P_{\alpha})^\vee$.
\end{proof}

\begin{definition}
We denote by $V(X)$ the scheme-theoretic support of  $(-1_{\widehat A} )^* R^g\widehat{\mathcal S}({\bf R}\Delta(\I_X(\Theta)))$, and call it the \emph{theta-dual} of $X$ (with respect to the principal polarization $\Theta$).
By Lemma \ref{identification}, set-theoretically we have
$$V( X) = \{ \alpha ~|~ h^0 (A, \I_X( \Theta) \otimes P_{\alpha}) \neq 0\} \subset \widehat A,$$
which via the principal polarization is identified with the locus $\{ a \in A ~|~ X \subset \Theta_a \}$
of theta-translates containing $X$. More precisely we have
\end{definition}

\begin{corollary}\label{line_bundle}
The sheaf  $(-1_{\widehat A} )^* R^g\widehat{\mathcal S}({\bf
R}\Delta(\I_X(\Theta)))$ is a line bundle on $V(X)$. In fact
$$(-1_{\widehat A} )^* R^g\widehat{\mathcal
S}({\bf R}\Delta(\I_X(\Theta)))\cong {\mathcal O}_{V(X)}(\Theta).$$
\end{corollary}
\begin{proof}
The first assertion follows from the previous Lemma: the fibers of
$(-1_{\widehat A} )^* R^g\widehat{\mathcal S}({\bf
R}\Delta(\I_X(\Theta)))$ have dimension either zero or one, and the
latter happens if and only if $\alpha \in V(X)$. For the last
assertion, note that we have a natural surjective homomorphism
$$R^g \widehat{\mathcal{S}} (\R \Delta (\OO_A(\Theta))) \longrightarrow R^g\widehat{\mathcal S}({\bf R}\Delta(\I_X(\Theta))).$$
Indeed, as in Lemma \ref{identification}, the fiber over $\alpha$ is given by the dual of the injection
$H^0 (A, \I_X(\Theta) \otimes P_{\alpha}) \hookrightarrow H^0 (A, \OO_A(\Theta) \otimes P_{\alpha})$.
Now $\R \Delta (\OO_A(\Theta)) \cong \OO_A(-\Theta)$ and
$(-1_{\widehat A} )^* R^g\widehat{\mathcal S}({\mathcal O}_A(-\Theta))\cong {\mathcal O}_{\widehat A}(\Theta)$, which gives the conclusion.
\end{proof}


Recall from \cite{ran} \S II that an equidimensional reduced subscheme $X$ is \emph{nondegenerate} in $A$ if the kernel of the restriction map
$H^0 (A, \Omega_A^d) \rightarrow H^0 (X_{reg}, \Omega_{X_{reg}}^d)$ is $0$, and \emph{geometrically nondegenerate} if the same kernel contains no nonzero decomposable $d$-forms. It is enough here to work with geometrically nondegenerate subschemes.\footnote{We thank O. Debarre for pointing this out.}  Given that the image of the difference map $\phi: X \times V(X) \rightarrow A$, $(x,y) \to x-y$, is contained in $\Theta$, the following is a special case of \cite{debarre2} Proposition 1.4.

\begin{lemma}\label{nondegenerate}
Let $X$ be a geometrically nondegenerate equidimensional reduced subscheme of $A$, of dimension $d$. Then
$${\rm dim}~V(X) \le g-d-1.$$
\end{lemma}

\begin{example}\label{w_d}
If $C$ is a smooth projective curve of genus $g$, consider the image $W_d$
of the $d$-th symmetric product of $C$ via some
Abel-Jacobi map. Then it is folklore, and straightforward to verify,
that set-theoretically $V(W_d) = - W_{g-d-1}$ up to translate.
We will see later that this actually holds scheme-theoretically,
and that a similar fact holds for the Fano surface (cf. \S8.1).
\end{example}

\section{Consequences of the generic vanishing criterion}

We consider now $X$ to be a closed reduced subscheme of $A$ of \emph{pure dimension} $d$.
We draw a number of sheaf-theoretic consequences
in the case when $\I_X(\Theta)$ is a $GV$-sheaf, based on the criterion in \S2. Recall that by Theorem \ref{gv}  this is equivalent to the fact that the complex $\R \Delta (\I_X(\Theta))$ satisfies $WIT_g$ which, together with Corollary \ref{line_bundle}, means that
\begin{equation}\label{interpretation}
{\bf R}\widehat{\mathcal{S}}( {\bf R}\Delta(\I_X(\Theta)))= R^g\widehat{\mathcal{S}}({\bf
R}\Delta(\I_X(\Theta))) [-g]  \cong (-1_{\widehat A})^* \OO_{V(X)}(\Theta) [-g].
\end{equation}

\begin{proposition}\label{dimension}
Assume that $\I_X(\Theta)$ is a $GV$-sheaf. Then:

\noindent
(a) $\dim V(X) \ge g-d-1$.

\noindent
(b) If equality is attained, i.e. $\dim V(X) =  g-d-1$, then $V(X)$ is Cohen-Macaulay and equidimensional if and only if $\I_X(\Theta)$ satisfies  $WIT_{d+1}$.
\end{proposition}
\begin{proof}
(a) Note first that
$${\rm codim}_{\widehat A} V(X)={\rm min}\{k\>|\> {\mathcal
E}xt_{\OO_{\widehat A}}^k(\OO_{V(X)},\OO_{\widehat A})\ne 0\}.$$
On the other hand, we claim that ${\bf R}\widehat{\mathcal{S}}(\I_X(\Theta)) \cong  {\bf R}\Delta(\OO_{V(X)}(\Theta))$, which implies for all $i$ that
\begin{equation}\label{wit}
{\mathcal E}xt_{\OO_{\widehat A}}^i(\OO_{V(X)}(\Theta),\OO_{\widehat A}) \cong
R^i\widehat{\mathcal{S}}(\I_X(\Theta)) .
\end{equation}
To prove the claim, start with ($\ref{interpretation}$) above and apply the functor $ (-1_{\widehat A})^*\circ  \R \Delta\circ [g]$ to both sides. On the right hand side we obtain ${\bf R}\Delta(\OO_{V(X)}(\Theta))$. On the left hand side we have
$$\R \Delta \big( (-1_{\widehat A})^* \R \widehat{\mathcal{S}} ~\R \Delta(\I_X(\Theta) [g])\big) \cong
\R \Delta \big(\R \Delta ~\R \widehat{\mathcal{S}} (\I_X(\Theta))\big ) \cong
\R \widehat{\mathcal{S}} (\I_X(\Theta)),$$ where for the first isomorphism we
used Grothendieck duality  (\ref{gd}).

Now since $X$ has dimension $d$, we have that $R^i \widehat{\mathcal{S}} (\OO_X(\Theta)) = 0$ for $i > d$, and so $R^i \widehat{\mathcal{S}} (\I_X(\Theta)) = 0$ for $i > d + 1$, by applying the Fourier-Mukai functor to
the exact sequence
$$0 \rightarrow \I_X(\Theta) \rightarrow \OO_A(\Theta) \rightarrow \OO_X(\Theta) \rightarrow 0.$$
By the above this implies that the dimension of $V(X)$ is at least $g-d-1$.

\noindent
(b) This follows directly from (\ref{wit}), as it is well known that the condition on
$V(X)$ is equivalent to the fact that
$${\mathcal E}xt_{\OO_{\widehat A}}^k(\OO_{V(X)},\OO_{\widehat A})\ne 0 \iff k = d+1.$$
\end{proof}

\begin{theorem}\label{V(X)}
Assume that  $\I_X(\Theta)$ is a $GV$-sheaf and $X$ is geometrically nondegenerate. Then:

\noindent
(a) $V(X)$ is Cohen-Macaulay of pure dimension $g-d-1$.

\noindent
(b) $\I_{V(X)} (\Theta)$ is a $GV$-sheaf.

\noindent
(c) $V(V(X)) = X$.

\noindent
(d) $X$ is Cohen-Macaulay.
\end{theorem}
\begin{proof}
(a) Note first that $\dim V(X) = g-d-1$, since by Proposition \ref{dimension}(a) it is at least $g-d-1$, while by Lemma \ref{nondegenerate} it is at most that much. But then by Proposition \ref{dimension}(b),
the conclusion holds if we show that $\I_X(\Theta)$ satisfies $WIT_{d+1}$.

This is the key point in the proof of the entire theorem: we apply Theorem \ref{gv} to the complex $\F={\bf R}\Delta(\I_X(\Theta))$ (so that ${\bf R}\Delta \F=\I_X(\Theta)$). We know that  $\F$ is $GV_{g-d-1}$ in a very strong sense: $R^i \widehat{\mathcal{S}} \F =0$ for
$i\ne g$, while ${\rm codim}_{\widehat A}{\rm Supp}( R^g \widehat{\mathcal{S}} \F) = {\rm codim}_{\widehat A} V(X) = d+1= g - (g-d-1)$. Thus in Theorem \ref{gv} we can take $k = d-g-1$. The implication (a) $\Rightarrow$ (c) in the theorem gives that $R^i \widehat{\mathcal{S}} (\I_X(\Theta)) = R^i \widehat{\mathcal{S}} (\R \Delta \F) = 0$ for $i < d+1$. On the other hand we noted above that we have
$R^i \widehat{\mathcal{S}} (\I_X(\Theta)) = 0$ for $i > d+1$, simply since $X$ is $d$-dimensional.

\noindent
(b) Again by Theorem \ref{gv}, this is the same as showing that ${\bf R}\Delta(\I_{V(X)}(\Theta))$ satisfies $WIT_g$. We prove more precisely that
\begin{equation}\label{old_dual}
\R \mathcal{S} ({\bf R}\Delta(\I_{V(X)}(\Theta))) \cong (-1_A)^* \OO_X(\Theta).
\end{equation}
Using Mukai's inversion theorem for the Fourier functor, i.e. equation ($\ref{inversion}$) in \S2, this statement is equivalent to showing that
\begin{equation}\label{old}
\R \widehat{\mathcal{S}} (\OO_X(\Theta)) \cong \R \Delta (\I_{V(X)} (\Theta)).
\end{equation}
Because of the short exact sequence
$0 \longrightarrow \I_X (\Theta) \longrightarrow \OO_A (\Theta) \longrightarrow
\OO_X (\Theta) \longrightarrow 0$
and its analogue for $V(X)$, this is in turn equivalent to
$$\R \widehat{\mathcal{S}} (\I_X(\Theta)) \cong \R \Delta (\OO_{V(X)} (\Theta)).$$
This was shown in the proof of Proposition \ref{dimension}.

\noindent
(c) The assertion follows immediately from (\ref{old_dual}) above, since as in Lemma
\ref{identification} the Fourier transform  of ${\bf R}\Delta(\I_{V(X)}(\Theta))$ is supported on $-V(V(X))$.

\noindent
(d) From (a) we know that $V(X)$ is equidimensional, and from (b) that $\I_{V(X)} (\Theta)$ is a $GV$-sheaf. Since we know already that  $\dim V(X) = g-d-1$ and $\dim V(V(X)) = \dim X = d$, we can deduce as in (a) (without assuming nondegeneracy for $V(X)$) that $\I_{V(X)} (\Theta)$ satisfies
$WIT_{g-d}$. Proposition \ref{dimension}(b) implies then that $X$ is Cohen-Macaulay.
\end{proof}

\begin{remark}
Since $X$ is Cohen-Macaulay, it is  a posteriori not hard to understand the main technical tool,
namely the object $\R \Delta (\I_X(\Theta))$, better: in fact it is represented by a
complex with two non-zero cohomology sheaves: $\OO_A(-\Theta)$ in
degree $0$ and $\omega_X(-\Theta)$ in degree $g-d-1$.
\end{remark}

\begin{corollary}\label{twice}
If $X$ is geometrically nondegenerate and $\I_X(\Theta)$ is a $GV$-sheaf, then $\omega_X(-\Theta)$ satisfies $WIT_d$ and
$$\widehat{\omega_X(-\Theta)} \cong (-1_{\widehat A})^* \I_{V(X)}(\Theta).$$
Also $\omega_{V(X)} (-\Theta)$ satisfies $WIT_{g-d-1}$ and
$$\widehat{\omega_{V(X)} (-\Theta)} \cong (-1_ A)^* \I_X (\Theta).$$
\end{corollary}
\begin{proof}
First of all,  since we've seen that $X$ and $V(X)$ are Cohen-Macaulay, they do have dualizing sheaves. The second assertion is the same as the first, with the roles of $X$ and $V(X)$ reversed. For the first, note that $\R \Delta (\omega_X (- \Theta)) \cong \OO_X(\Theta) [g-d]$. Hence the statement follows from the isomorphism (\ref{old}) above and Grothendieck duality (\ref{gd}), as in the proof of Proposition \ref{dimension}.
\end{proof}

\begin{remark}\label{symmetry}
The results of this section imply that if the $GV$ condition is satisfied for $\I_X(\Theta)$ then the role of the schemes $X$ and $V(X)$ is symmetric: any type of result we are interested in which is true for $X$ holds also for $V(X)$, and conversely. Note that Theorem \ref{bis} below implies that $V(X)$ is in fact \emph{nondegenerate}, and so the same thing holds for $X$.
\end{remark}

\section{Minimal cohomology classes}

In this section we prove a result which, by the symmetry
noted in Remark \ref{symmetry}, is equivalent to Theorem \ref{b}.
Besides the results of the previous section, we will use applications of
the Grothendieck-Riemann-Roch theorem.

\begin{theorem}\label{bis}
Let $X$ be a geometrically nondegenerate closed reduced subscheme of pure dimension $d$ of a ppav $(A,\Theta)$ of dimension $g$, such that $\I_X(\Theta)$ is a $GV$-sheaf.
Then $[V(X)] = \frac{\theta^{d+1}}{(d+1)!}$.
\end{theorem}
\begin{proof}
Recall from Corollary \ref{twice} that $(-1_A)^* \omega_X (-\Theta)$ satisfies $WIT_d$,
and its Fourier-Mukai transform is $\I_{V(X)} (\Theta)$, which satisfies $WIT_{g-d}$. We will use the relationship in cohomology between the Chern character of a $WIT$-sheaf and that of its Fourier-Mukai transform, established by Mukai \cite{mukai2} Corollary 1.18\footnote{This says that if $\F$ is a sheaf on $A$ satisfying $WIT_j$, then ${\rm ch}_i (\widehat \F) = (-1)^{i + j} PD_{2g-2i} ({\rm ch}_{g-i} (\F))$.} as an application of
Grothendieck-Riemann-Roch. In our situation, for the $(d+ 1)$-st component of the Chern character (indexed by the codimension), this is written as
$${\rm ch}_{d+1} (\I_{V(X)} (\Theta)) = (-1)^{g+1} PD ({\rm ch}_{g-d-1} ((-1_A)^* \omega_X(-\Theta))),$$
where $PD$ denotes the Poincar\' e duality isomorphism, and we ignore notation for the inclusion map of $X$ into $A$.

We now compute the two sides. On one hand, since $\omega_X(-\Theta)$ is supported on $X$, which
has dimension $d$, it is standard that  ${\rm ch}_{g-d-1} ((-1_A)^* \omega_X(-\Theta)) = 0$ since
the Chern character ${\rm ch} (\omega_X (-\Theta))$ has no non-zero components in degree less than the codimension $g-d$. (For example this is a simple consequence of Grothendieck-Riemann-Roch -- cf. e.g. \cite{fulton} Example 18.3.11.)

On the other hand
$${\rm ch} (\I_{V(X)}(\Theta)) = {\rm ch} (\OO_{\widehat A} (\Theta)) - {\rm ch} (\OO_{V(X)} (\Theta)) = $$
$$= \big(1 + \theta + \frac{\theta^2}{2!} + \ldots \big) - \big([V(X)] + \ldots \big),$$
The second term starts in $CH^{d+1}$; we apply the same reasoning as in the previous paragraph, in the more precise form saying that ${\rm ch}(\OO_{V(X)})$ starts with the class
$[V(X)]$ (cf. \cite{fulton} Example 15.2.16), using the fact that by Theorem \ref{V(X)}(a), $V(X)$ is of pure dimension $g-d-1$. We obtain that ${\rm ch}_{d+1} (\I_{V(X)} (\Theta)) = \frac{\theta^{d+1}}{(d+1)!} - [V(X)]$. This has to be $0$ by the paragraph above, which gives the conclusion.
\end{proof}

\section{Characterization of Jacobians via generic vanishing for extremal dimensions}

In this section we prove Theorem \ref{c}. At this stage it is in fact an immediate corollary of
Theorem \ref{bis} and the Matsusaka-Hoyt criterion \cite{hoyt}, Theorem 1.

\begin{proof}(\emph{of Theorem \ref{c}.})
By Theorem \ref{bis} and Theorem \ref{V(X)} we know that both $X$ and $V(X)$ have minimal class.
If $X$ is  a curve, then the statement follows directly from the Matsusaka criterion.

If $X$ has codimension $2$, then $V(X)$ is a $1$-cycle of minimal class.
By the Matsusaka-Hoyt criterion we obtain that $V(X)$ is an irreducible Abel-Jacobi embedded curve $C$, since otherwise $(A, \Theta)$ would have to split into a product of polarized Jacobians, which contradicts indecomposability. Now given any Abel-Jacobi embedded curve $C$, by Example \ref{w_d} we have $V(C) = - W_{g-2} (C)$ in the dual Jacobian. Since $V(V(X)) = X$ by Theorem
\ref{V(X)}(c), this concludes the proof. (The entire discussion is of course up to $+$ or $-$ a translate).
\end{proof}

\begin{remark}
The exact same proof tells us what happens in case $(A, \Theta)$ is decomposable as well. The
Matsusaka-Hoyt theorem implies that $A$ is a product of Jacobians and either $X$ or $V(X)$, according to which has dimension $1$, is a reduced $1$-cycle which projects to an Abel-Jacobi embedded curve on each one of the factors.
\end{remark}

\section{Complements}

\subsection{The theta-dual of $W_d$ and of the Fano surface}
Going back to Example \ref{w_d}, we can now say that $V(W_d)=-W_{g-d-1}$ \emph{scheme-theoretically}. The equality is known set-theoretically, and by Theorem \ref{V(X)}(a), $V(W_d)$ has no embedded components. Since $[W_{g-d-1}]= \frac{\theta^{d+1}}{(d+1)!}$, the conclusion follows from Theorem \ref{b}.

It was shown in \cite{cg} that for the Fano surface of lines $F$ in the intermediate Jacobian of a smooth
cubic threefold we have $F - F = \Theta$ set-theoretically, and $[F] = \frac{\theta^3}{3!}$. The first equality gives $F \subset V(F)$, so given the results above $V(F)$ has dimension $2$ and $F$ as a component. A priori it could have other components too. The results of \cite{hoering} imply however that $\I_F(\Theta)$ is $GV$, which again by Theorem \ref{V(X)}(a) and Theorem \ref{b} gives that $V(F)$ is equidimensional and $[V(F)] = \frac{\theta^3}{3!}$.  This implies that $V(F) = F$.

\subsection{A remark on the Beauville-Debarre-Ran question}
In \S6 we saw that (2) implies (1) in Conjecture \ref{a}. The full
equivalence of (1) and (2) would have an amusing consequence for the equivalence of (1) 
and (5): it would imply that it is enough to check it only for subvarieties of
dimension $d \le [\frac{g}{2}]$ (or only for subvarieties of
dimension $d \ge [\frac{g}{2}]$), i.e. essentially for half the
dimensions involved. This follows from the proof of Theorem \ref{b},
where we saw that $\I_X(\Theta)$ is a $GV$-sheaf if and only if
$\I_{V(X)} (\Theta)$ is a $GV$-sheaf, combined with the fact that on
a Jacobian $J(C)$ we have $V(W_d (C)) = -W_{g-d-1} (C)$ up to
translate (cf. Example \ref{w_d}).

\subsection{An addition to Conjecture \ref{a} involving the theta-dual $V(X)$}
Finally we make a comment, independent of the arguments above, emphasizing the fact that one can try to detect subvarieties of minimal cohomology class based directly on properties of $V(X)$.
If we again identify $V(X)$ with the locus in $A$ parametrizing theta-translates containing $X$, then
we have a difference map
$$\phi: X \times V(X) \rightarrow X - V(X) \subset \Theta \subset A.$$
Consider now the equalities
\begin{equation}\label{pontrjagin}
{\rm deg}(\phi)\cdot[X - V(X)] = [X] *  [V(X)] = {{g-1}\choose {d}} \theta.
\end{equation}
where $*$ denotes the Pontrjagin product on cohomology. These are satisfied for the difference maps $W_d \times W_{g-d-1} \rightarrow W_d - W_{g-d-1}$ and $F \times F \rightarrow F-F$, as for any $d$ we have the formula
$$\frac{\theta^{g-d}}{(g-d)!} * \frac{\theta^{d+1}}{(d+1)!}  =  {{g-1}\choose {d}} \theta$$
(cf. \cite{lb} Corollary 16.5.8). On the other hand, if (\ref{pontrjagin}) holds for some subscheme $X$ of dimension $d$, then $X$ has minimal class if and only if $V(X)$ has minimal class.\footnote{Here is a brief explanation: via the Fourier transform $F$ on cohomology, which is an isomorphism, the Pontrjagin product corresponds to the cup product. On the other hand, $F (\frac{\theta^a}{a!}) = \frac{\theta^{g-a}}{(g-a)!}$. But taking cup-product with this class is an isomorphism by the Lefschetz theorem. For details on all of this cf. e.g. \cite{beauville2}, or \cite{lb} \S16.3.}
This suggests adding a more geometric conjecture to the list in Conjecture \ref{a} -- namely  under the hypotheses of Conjecture \ref{a}, the following should be equivalent, and equivalent  to the statements there:
\begin{enumerate}
\item ${\rm dim}~V(X) = g-d-1$.
\item $V(X)$ has minimal cohomology class $[V(X)] = \frac{\theta^{d+1}}{(d+1)!}$.
\item $[X] *  [V(X)] = {{g-1}\choose {d}} \theta$.
\end{enumerate}

Note that (3) and (2) obviously imply (1). The key point is (1): subvarieties of
minimal class should be characterized among nondegenerate subvarieties by the fact that
 their theta-dual variety has \emph{maximal dimension}.

\providecommand{\bysame}{\leavevmode\hbox to3em{\hrulefill}\thinspace}

\end{document}